\documentclass{proc-l}

 \usepackage{amssymb, epsfig, amssymb, latexsym}

\newtheorem{theorem}{Theorem}[section]
\newtheorem{lemma}[theorem]{Lemma}

\theoremstyle{definition}

\theoremstyle{remark}
\newtheorem{remark}[theorem]{Remark}

\numberwithin{equation}{section}

\newcommand{\p}{\partial}
\newcommand{\D}{\Delta}
\newcommand{\om}{\omega}
\newcommand{\Om}{\Omega}
\renewcommand{\phi}{\varphi}
\newcommand{\e}{\epsilon}
\renewcommand{\a}{\alpha}

   \newcommand{\EX}{{\Bbb{E}}}
   \newcommand{\PX}{{\Bbb{P}}}

\begin{document}

\title[Stochastic parameterization]{Stochastic parameterization  for large eddy simulation of geophysical flows    }


\author{Jinqiao Duan}
\address{Department of Applied Mathematics\\ Illinois Institute of Technology \\
  Chicago, IL 60616, USA.}
\curraddr{} \email{duan@iit.edu}
\thanks{Partly supported by the NSF Grants DMS-0209326  \& DMS-0542450.}

\author{Balasubramanya T. Nadiga}
\address{Los Alamos National Laboratory\\
 MS-B296, Los Alamos, NM 87545, USA.}
\curraddr{} \email{balu@lanl.gov}
\thanks{Partly supported  by the LDRD program at Los Alamos National Laboratory.}


\date{October 30, 2005 and December 30, 2005 (Revised version)}

\dedicatory{}

\commby{Ed Waymire}

\begin{abstract}
  Recently, stochastic, as opposed to deterministic, parameterizations
  are being investigated to model the effects of unresolved subgrid
  scales (SGS) in large eddy simulations (LES) of geophysical flows.
  We analyse such a stochastic approach in the barotropic vorticity
  equation to show that (i) if the stochastic parameterization
  approximates the actual SGS stresses, then the solution of the stochastic
  LES approximates the ``true" solution at appropriate scale sizes; and that
  (ii) when the filter scale size approaches zero, the solution of the
  stochastic LES approaches the true solution.

\end{abstract}

\maketitle


\section{Motivation}

The immense number of degrees of freedom in large scale turbulent
flows as encountered in the world oceans and atmosphere makes it
impossible to simulate these flows in all their detail in the
foreseeable future. On the other hand, it is essential to
represent these flows reasonably accurately in Ocean and
Atmospheric General Circulation Models (OGCMs and AGCMs) so as to
improve the confidence in these model components of the earth
system in ongoing effort to study climate and its variability
(e.g., see \cite{ipcc}). Furthermore, it is very often the case
that in highly resolved computations, a rather disproportionately
large fraction of the computational effort is expended on the
small scales (e.g., see \cite{pope}) whereas a large fraction of
the energy resides in the large scales (e.g., see \cite{nastrom}).
It is for these reasons that the ideas of Large Eddy Simulation
(LES)---wherein the large scale unsteady motions driven by
specifics of the flow are explicitly computed, but the small (and
presumably more universal \cite{MK00}) scales are modelled---are
natural in this context.

Given our interest in large scale geophysical flows with its small
vertical to horizontal aspect ratio, we restrict ourselves to
two-dimensional or quasi two-dimensional flows.  Previous models
of the small scales in how they affect the large scales in the
momentum equations or equivalently the vorticity equation in
incompressible settings have mostly been confined to an enhanced
eddy viscosity or nonlinear eddy viscosity like that of
Smagorinsky or biharmonic viscosity (e.g., see \cite{smag1,smag2,
leith1,leith2, holland,bleck}). Given the non-unique nature of the
small scales with respect to the large scales \cite{pope}, the
aforementioned use of deterministic and dissipative closures seem
rather highly restrictive. On the other hand, it would seem
desirable to actually represent a population of eddies that
satisfy overall constraints of the flow rather than make flow
specific parametric assumptions. This has led to recent
investigations of the possibility of using stochastic processes to
model the effects of unresolved scales in geophysical flows (e.g.,
see \cite{majda1, majda2}).  More recently, subgrid scale (SGS)
stresses have been analysed in simple but resolved flows as a
possible way to suggest stochastic parameterizations as e.g, in
\cite{Berloff,Berloff2,Nadiga1, Nadiga2}.  These efforts have been
preceded, of course, by various attempts to model anomalies in
geophysical flow systems as linear Langevin equations (e.g,
\cite{Hasselmann, Arn00, herring, DelSole-Farrell}) and the
analysis of stochastic models in isotropic and homogeneous three
dimensional turbulence (e.g., \cite{kraichnan,Monin}).

In this paper, we analyze the stochastic approach to parameterization
in the barotropic vorticity equation and show that (i) if the
stochastic parameterization approximates the SGS stresses, then the
stochastic large eddy solution approximates the ``true" solution at
appropriate scale sizes; and that (ii) when the filter scale size
approaches zero, then the solution of the stochastic LES approaches
the true solution.  In the next section we present a set of
computations that demonstrates the use of stochastic parameterizations
and in \S 3, we prove the main results on approximation and
convergence of LES solutions using stochastic parameterizations.

\section{Stochastic parameterization}
We consider the simple setting of the beta-plane barotropic
vorticity equation (equivalently the two-dimensional (2D)
quasi-geostrophic (QG) model) \cite{Gill, Pedlosky}:
\begin{eqnarray} \label{QG}
 q_t + J(\psi, q) + \beta \psi_x = f(x,y,t) + \nu \Delta q - r q,
\end{eqnarray}
on a bounded domain $D$ with piecewise smooth boundary $\partial D$.
Here the vorticity $q (x,y,t)$ is given in terms of streamfunction
$\psi (x,y,t)$ by $q=\Delta \psi$. $\beta$ is the meridional gradient
of the Coriolis parameter, $\nu$ $>$ $0$ the viscous dissipation
constant, $r$ $>$ $0$ the Ekman dissipation constant and $f(x,y,t)$
the wind forcing.  The forcing $f$ is always assumed to be mean-square
integrable both in time and in space.  In addition, $\Delta$ $=$
$\p_{xx}+\p_{yy} $ is the Laplacian operator in the plane and $J(h,g)$
$=$ $h_xg_y -h_yg_x$ is the Jacobian operator. The boundary condition
(BC) is $ q =0, \; \psi =0\;$ on $\partial D$ and initial condition
(IC) is $ q (x,y,0)= q_0(x,y). $

Fine mesh simulations ($q$) are used to obtain the benchmark solution
$ \bar{q} $ through convolution with a spatial filter
$G_{\delta}(x,y)$, with spatial scale $\delta>0$:
\begin{eqnarray} \label{benchmark}
\bar{q} (x,y,t) : = q* G_{\delta}
\end{eqnarray}
We use a Gaussian filter
\cite{Berselli},
$
G_{\delta}(x,y)= \frac{1}{\pi\delta^2} e^{-\frac{
    (x^2+y^2)}{\delta^2}}, $ where $\delta>0$ is the filter size and
the filter is such that (1) $q* G_{\delta}$ is infinitely
differentiable in space and (2) $q*G_{\delta} \to q$ as $\delta \to
0$ in $L^2(D)$.  Note that the Fourier transform of $G_{\delta}$ is
$\widehat{G}_{\delta}(k_1,k_2)= e^{-\frac{\delta^2 (k_1^2+k_2^2)}{4 }}, $
and that $ \widehat{q*G}_{\delta} (k_1, k_2,
t)=\widehat{G}_{\delta}(k_1,k_2) \widehat{q} (k_1, k_2, t).  $

On convolving (\ref{QG}) with $G_{\delta}$
the large eddy solution $\bar{q}$ is seen to satisfy
\begin{eqnarray} \label{bar}
 \bar{q}_t + J(\bar{\psi},\bar{q}) + \beta \bar{\psi}_x=
 \bar{f}(x,y,t) + \nu \Delta \bar{q} - r \bar{q}  + R(q),
\end{eqnarray}
where  the SGS stress term $R(q)$ is defined as
\begin{eqnarray} \label{EF}
 R(q):=J(\bar{\psi},\bar{q})- \overline{J(\psi, q)}.
\end{eqnarray}
Note that since $\bar{\bar{q}} \ne \bar{q}$, the SGS stress term
$R(q)$ above is more than the Reynolds stress and the SGS stress is
usually further divided into three components, the explicit Leonard
stress, and the cross stress and the SGS Reynolds stress that require
further modeling. However, for our purposes, we will consider $R(q)$
in its entirety.  Since $R(q)$ depends on $q$ as well as $\bar{q}$,
the equation (\ref{bar}) is not a closed system. We need to model or
prescribe $R(q)$ in terms of resolved quantity $\bar{q}$.  On the
other hand, $R(q)$ may be explicitly diagnosed from a fully-resolved
run, given the filter.

An analysis of $R(q)$ reveals that its time-mean is much smaller
than its standard deviation, and that its temporal behavior is
highly irregular, leading to the possibility of approximating
$R(q)$ by a suitable stochastic process $\sigma(\bar{q},\om)$
(defined on a probability space $(\Om, \mathcal{F}, \mathbb{P})$,
with
 $\om \in \Om$, the sample space, $\sigma-$field $\mathcal{F}$ and probability
measure $ \mathbb{P} $).
With such a putative stochastic closure, the LES model becomes a random partial
differential equation (PDE) for $Q \sim \bar{q}$:
\begin{eqnarray} \label{LES}
 Q_t + J(\Psi,Q) +\beta \Psi_x= \bar{f}(x,y,t) + \nu \Delta Q - r Q + \sigma
 (Q, \omega)  ,
\end{eqnarray}
where $Q=\Delta \Psi$ and $\bar{f}(x,y,t) : = f* G_{\delta}$, with
BC: $
 Q =0, \; \Psi =0\;$    on    $\partial D$
 and IC:
$ Q (x,y,0)= \bar{q}_0(x,y). $ Note that when the stochastic
process $\sigma(Q,\omega)$ does not depend explicitly on $Q$---as
is the case when for example, either $R(q)$ or some of its
statistical properties are used to construct $\sigma(\omega)$---we
end up with an additive stochastic closure. The more general case
of the stochastic closure wherein there is an explicit dependence
on the state of the system $Q$, corresponds to a multiplicative
stochastic closure.

While one can get an idea of the stochastic forcing to be used to
represent the effects of unresolved subgrid scales in the LES runs by
analysing resolved runs at the scale of the LES computation, the
selection of a specific functional form for the stochastic
parameterization is beyond the scope of the present article. We aim at
providing a quantitative guide to selecting the stochastic
parameterization.

\begin{figure}[h]
\includegraphics[width=\textwidth]{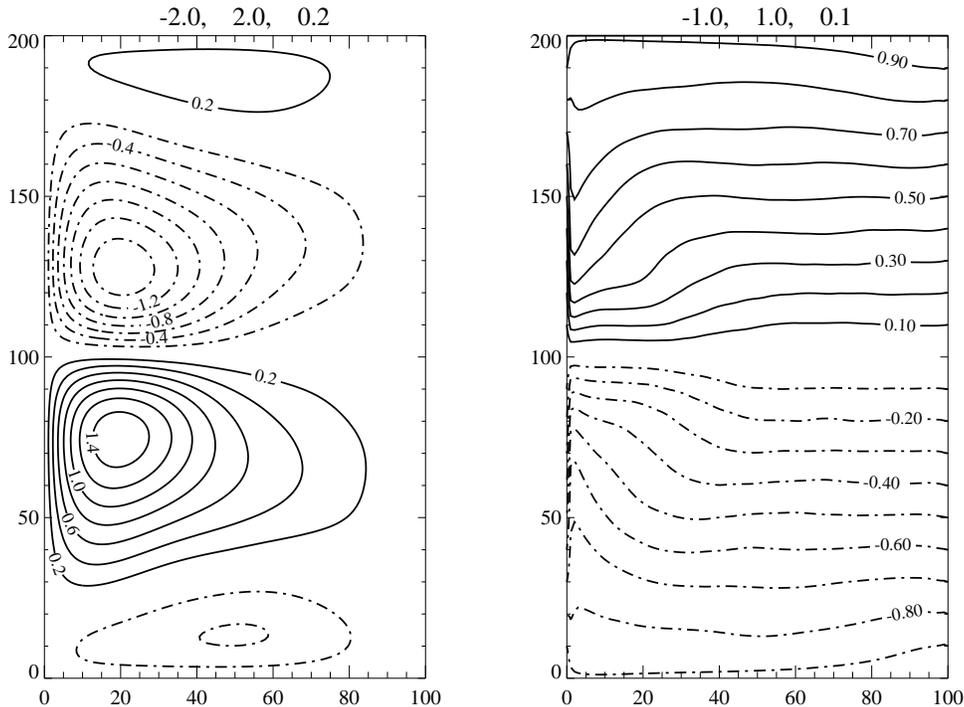}
\caption{\label{resolved}Resolved simulation. Contour
plots of the average over the attractor (time-average) of
streamfunction and potential-vorticity.
Contour intervals are the same in all the figures and are
indicated on top of the figures in the form (min, max, increment).
}
\end{figure}

\subsection{Numerical experiments}
We now briefly present a set of computations using the beta-plane
barotropic vorticity equation (\ref{QG}) in a rectangular
midlatitude basin. Finite differencing in space is used along with
Runge-Kutta time stepping, and other details of the setup may
be found elsewhere \cite{Nadiga3}. The steady forcing is uniform
in the zonal ($x$) directions and sinusoidal in the meridional
($y$) direction corresponding to a double-gyre wind forcing. At
the parameter values that we are presently consider, the
circulation is highly variable, but statistically stationary. We
therefore consider long time averages over the attractor in place
of ensemble averaging (over $\omega$).

Fig.~1 shows the contour plots of the time average of
streamfunction and potential vorticity as they emerge in the
resolved computations. A discussion of the phenomenology of this
circulation may be found in \cite{Nadiga3}. This simulation is
then analysed using a Gaussian filter with a width that is four
times the grid spacing of the resolved computations to obtain the
SGS stress $R(q)$.

\begin{figure}
\includegraphics[width=\textwidth]{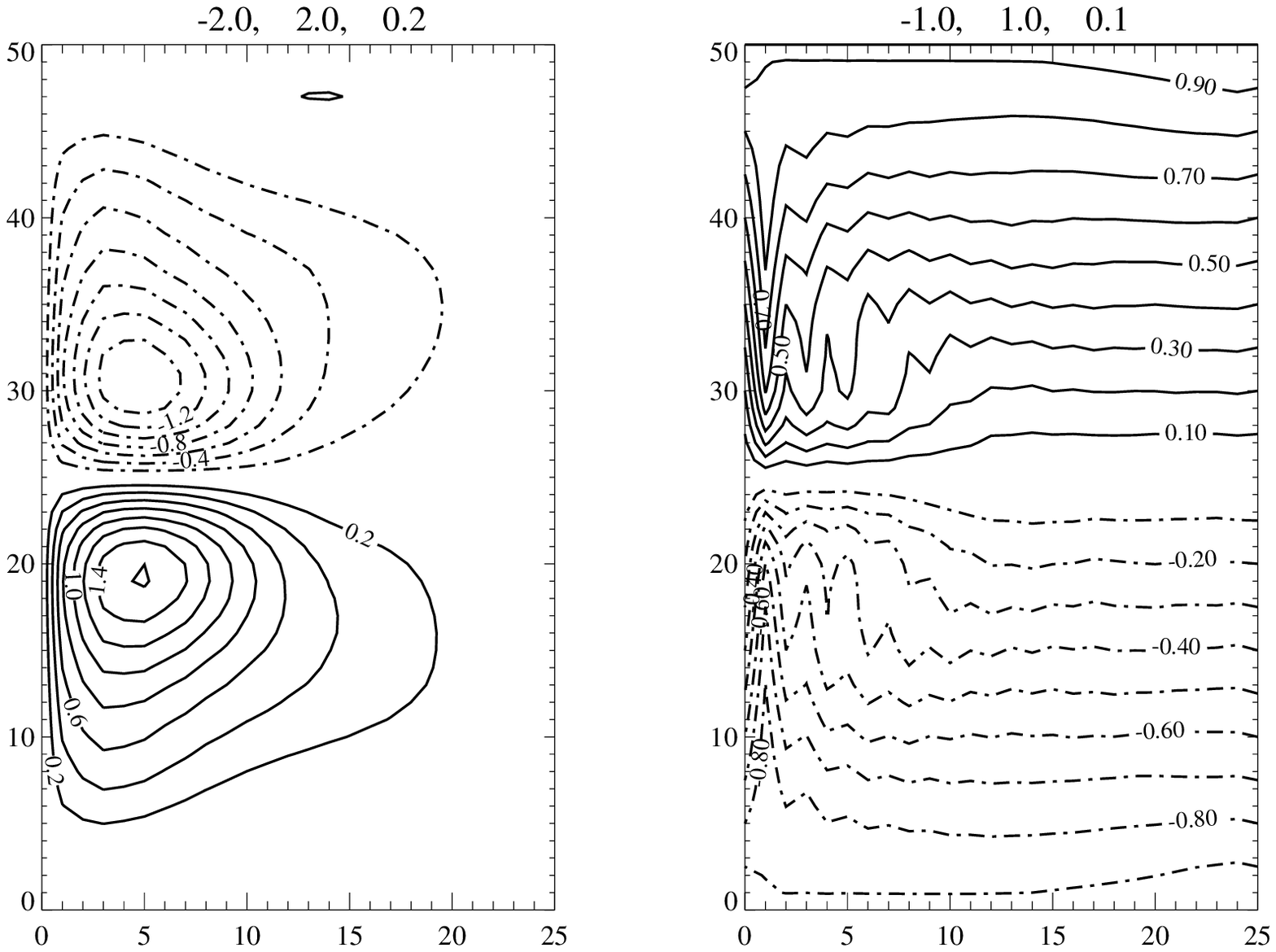}
\caption{\label{unresolved}Coarse unresolved simulation. Contour
plots of streamfunction and potential-vorticity, as in Fig.~1. The stochastic
parameterization term $\sigma(Q,\omega)$ is set to 0.}
\end{figure}

Next we consider a pair of coarse-scale simulations in which the
grid spacing is four times that of the resolved computation in
both directions. Fig.~2 shows the time average of the
streamfunction and potential vorticity as emerges from the
coarse-scale computation when $\sigma(Q,\omega)$ in (\ref{LES}) is
set to zero.  The main differences, as compared to the resolved
runs, clearly are the absence of the outer gyres in the
streamfunction field and the large amplitude grid scale
oscillation in the potential vorticity field.

In the next case, we set the statistics of $\sigma(Q,\omega)$
(viz., its spatial and temporal correlation functions, amplitude
and probability distribution function) identical to those of
$R(q)$ previously diagnosed from the resolved simulation. This we
do by using the actual time history of the SGS stress $R(q)$ in a
coarse-resolution run in which the initial condition is slightly
perturbed from the initial condition of the resolved run from
which $R(q)$ is diagnosed. Given the highly chaotic nature of the
flow, the effect of the initial perturbation is to quickly lead to
a complete decorrelation of the SGS stress forcing term $R(q)$
from the state of the system $Q$.  Thus, the actual time history
of the diagnosed SGS stress $R(q)$ supplied to the
coarse-resolution run acts effectively as an additive stochastic
closure $\sigma(\omega)$ of the SGS stresses in this LES. The time
averages of the circulation for this case is shown in Fig.~3.
Comparing Fig.~2 and Fig.~1 and Fig.~3 and Fig.~1, it is clear
that the differences between the resolved case and the case with
$\sigma(Q,\omega)$ similar to $R(q)$ is much smaller than the
differences between the resolved case and the case with
$\sigma(Q,\omega)=0$.

\begin{figure}
\includegraphics[width=\textwidth]{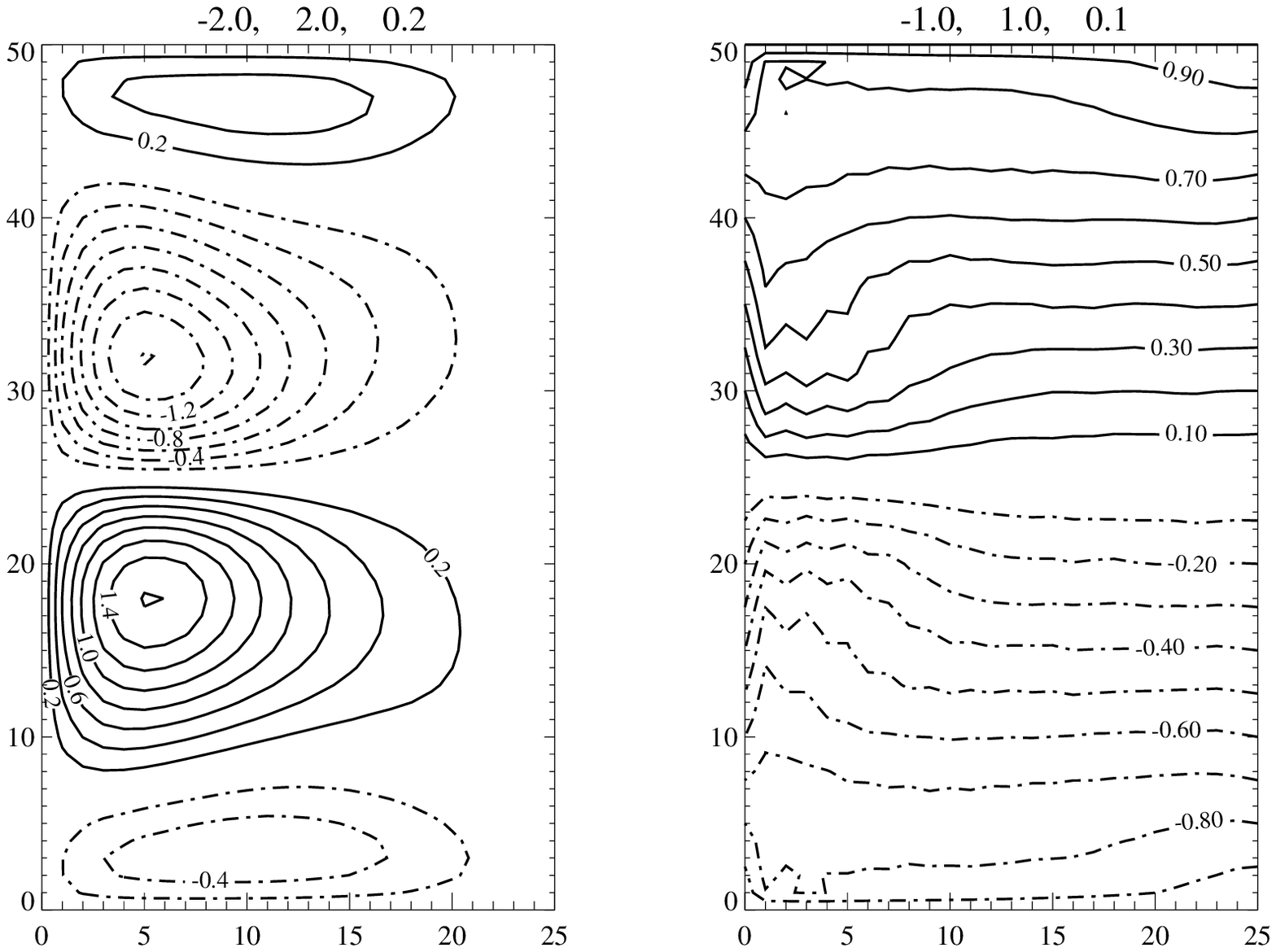}
\caption{\label{stochastic}LES with stochastic
parameterization. Contour plots of  streamfunction and potential-vorticity, as
in Fig.~1. The statistics of $\sigma(Q,\omega)$ are identical to those of $R(q)$.}
\end{figure}

These numerical experiments suggest that stochastic parameterizations
can provide good representations of subgrid scales. So, the question
then is as to how such an approach can be justified.

In an attempt to answer this question, we will prove, following the
approaches in \cite{Berselli,Kaya} in our stochastic context, that,
under appropriate conditions on the stochastic parameterization
(that appear easier to check in our case, c.f., Theorem 1 part
(i) below and the Assumption A2 in \cite{Kaya}):
\begin{eqnarray}
 \EX \| \bar{q} - Q\|^2  & \leq  & C(\nu, r,  q_0, T) \cdot \; \EX
 \int_0^T
 \|R(q)-\sigma(Q, \om)\|^2, \;\; 0\leq t \leq T,  \\
\EX \|q-Q\|^2  & \to &  0, \; \mbox{as} \;\; \delta \to 0,\;\;
0\leq t \leq T.
 \end{eqnarray}
 where $C(\cdot)>0$ is a constant, $[0, T]$ is the computational time
 interval, $ \EX (Z(\om)):= \int_{\Om} Z(\om) d\PX(\om)$, and $\|
 \cdot\|$ is the (spatial) norm in the space $L^2(D)$ of spatially
 mean-square integrable functions: $ L^2(D):= \{f:\;\; \|f\|
 =\sqrt{\int_D f(x,y)dxdy} < \infty \}. $


\section{Main results}

Standard abbreviations $L^2$ $=$ $L^2(D)$, $H^k_0$ $=$ $H^k_0(D)$,
$k = 1, 2,\ldots$,  are used for the common Sobolev spaces in
fluid mechanics \cite{Temam, Constantin}, with $<\cdot, \cdot>$
and $\| \cdot \|$ denoting the  usual (spatial) scalar product and
norm, respectively, in $L^2(D)$:
$$
<f, g>:=\int_D fg dxdy, \;\;\;\;\;\;  \|f\|:=\sqrt{<f, f>}
=\sqrt{\int_D f(x,y)dxdy}.
$$

 We need the following properties and
estimates (see also \cite{Dymnikov,Duan-Kloeden}) of the Jacobian
operator $J: H_0^1 \times H_0^1 \rightarrow  L^1$:
\begin{eqnarray*}
\int_D J(f,g) h \,dxdy   =   - \int_D J(f,h) g \,dxdy,
\int_D J(f,g) g \,dxdy   =   0,\quad\hbox{and}\\
    \label{estimate0}
\left|\int_D J(f,g)  \,dxdy \right| \leq  \|\nabla f\|
\|\nabla g\| \label{estimate2}
\quad\hbox{for all $f$, $g$, $h$ $\in$ $H^1_0$,}\\
\left|\int_D J(\Delta f, g) \Delta h \,dxdy \right|
 \leq \sqrt{\frac{2|D|}{\pi}} \|\Delta f\|  \; \|\Delta g\|  \;\|\Delta h\|
\label{estimate2.5}
\quad\hbox{for all $f$, $g$, $h$ $\in$ $H^2_0$}.
\end{eqnarray*}
In addition, we recall the Poincar\'e,  Young, and Gronwall inequalities below:
\begin{eqnarray*}
\|g\|^2 = \int_D g^2(x,y) \,dxdy \leq  \frac{|D|}{\pi} \int_D
|\nabla g|^2 \,dxdy = \frac{|D|}{\pi} \|\nabla g\|^2
\quad\hbox{for $g$ $\in$ $H^1_0$}.\\
AB \leq \frac{\e}2 A^2 + \frac{1}{2 \e}B^2
\quad\hbox{for $A,B \geq 0$ and $\e >0$}.
\end{eqnarray*}
Assuming that $y(t)\geq 0$, $g(t)$ and
$h(t)$ are integrable, if $\frac{dy}{dt}\leq g y +h$ for $t \geq
t_0$, then
\begin{equation*}
y(t)\leq y(t_0)e^{\int_{t_0}^t g(\tau)d\tau} +\int_{t_0}^t h(s)
e^{-\int_{t}^s g(\tau)d\tau} ds, \;\; t \geq t_0.
\end{equation*}

\begin{lemma} (\cite{Duan-Kloeden}) \label{lemma} The   quasi-geostrophic motion described by (\ref{QG})
satisfies the enstrophy estimate :
\begin{eqnarray} \label{estimate}
   \| q\|^2  & \leq  & \|q_0\|^2 e^{-2\a t}
   +\frac1{r} \int_0^t \|f(s)\|^2 e^{2\a(s-t)}ds ,  \;\;0 \leq t <\infty\\
 & \leq  & \|q_0\|^2 e^{ 2|\a|T}
  +\frac1{r} e^{ 2|\a|T} \int_0^T \|f(s)\|^2 ds ,  \;\;  0 \leq t \leq T
 \end{eqnarray}
where $\a =\frac{r}2 + \frac{\pi \nu}{|D|} - \frac12 |\beta|
\left(\frac{|D|}{\pi} +1\right)$  with $|D|$ denoting the area of
the domain $D$. Note that $\a$ is   positive in the case of no
rotation ($\beta=0$).
\end{lemma}

\begin{theorem} (i) \textbf{Stochastic Approximation}:

 If the stochastic paramerization  $\sigma(Q,\om)$ is such that
\begin{eqnarray} \label{condition}
  \int_0^T \|\sigma(Q,\om)\|^2dt \leq M(T),  \;\;
  \mbox{  almost surely for}\; \om \in \Om,
 \end{eqnarray}
 for some  constant $M>0$ depending on computational time interval,
 then
\begin{eqnarray}
 \EX \| \bar{q} - Q\|^2  & \leq  & C(\nu, r,  q_0, T) \cdot \; \EX
 \int_0^T\|R(q)-\sigma(Q, \om)\|^2 dt,  \;\;  0 \leq t \leq T,
 \end{eqnarray}
 for any fixed time interval $0\leq t \leq T$.

 This implies that, if the stochastic paramerization  $\sigma(Q,\om)$
 approximates the SGS stress $R(q)$, then the LES solution $Q$ approximates $\bar{q}$, in mean-square sense.
\\
 (ii) \textbf{Scale convergence}:
 If the stochastic paramerization  $\sigma(Q,\om)$ satisfies
\begin{eqnarray} \label{variance}
\EX \int_0^T \|\sigma (Q, \om)\|^2 dt \to 0,\;\; \; \mbox{as}\;\;
\delta \to 0,
\end{eqnarray}
 for all LES solutions $Q$  of \ref{LES}), then
\begin{eqnarray}
\EX \|q-Q\|^2  & \to &  0, \;\;\; \mbox{as} \;\;\; \delta \to 0,\;
\;\; 0<t< T.
 \end{eqnarray}
This implies that, if the stochastic paramerization
$\sigma(Q,\om)$ becomes smaller (collectively in
computational time interval) as the cut-off scale size $\delta$
decreases, then the LES solution $Q$ approximates the
original solution $q$ better, in mean-square sense.
\end{theorem}
\begin{remark}
Condition (\ref{condition}) means that the
stochastic paramerization $\sigma(Q,\om)$ is square-integrable in
time and space, and
 its norm in the space $L^2((0, T); L^2(D))$ is almost surely bounded on the computational
 interval.
\end{remark}
\begin{remark}
Condition (\ref{variance}) means that the variance of the
stochastic paramerization $  \sigma (q,\om)$,  collectively in the
finite time interval of numerical simulation, becomes smaller and
smaller as the cut-off scale size $\delta$ decreases.
\end{remark}
To prove part (i), denote $U= \bar{q} -Q$, so that $U= \D (\bar{\psi}
-\Psi)$. Note that $U(0)=0$. Subtracting (\ref{bar}) from
(\ref{LES}), we see that $U$ satisfies
\begin{equation} \label{U}
 U_t =- J(\bar{\psi},\bar{q}) +  J(\Psi,Q) -\beta
 (\bar{\psi}_x-\Psi_x)
  + \nu \Delta U - r U +  [R(q)- \sigma(Q, \omega)].
\end{equation}
Multiplying this equation by $U$,   integrating over $D$ and
noting that $\bar{q}=U+Q$, we obtain
\begin{eqnarray} \label{Udot}
 \frac12 \frac{d}{dt}\|U\|^2 =- \int_D J(\bar{\psi}-\Psi,Q)U
 -\beta\int_D  (\bar{\psi}_x-\Psi_x)U \\
  - \nu\|\nabla U\|^2 - r \|U\|^2 +  \int_D [R(q)- \sigma(Q, \omega)]U.
\end{eqnarray}

Note that, for $0 \leq t \leq T$,
\begin{eqnarray}
& & \left|\int_D J(\bar{\psi}-\Psi,Q)U  \,dxdy \right|  \nonumber\\
 & \leq & \sqrt{\frac{2|D|}{\pi}} \|Q\|  \; \|U\|^2  \nonumber \\
& \leq & \sqrt{\frac{2|D|}{\pi}} \; \{\|\bar{q}_0\|^2 e^{ 2|\a|T}
  +\frac1{r} e^{ 2|\a|T} \int_0^T (\|\bar{f}(s)\|^2 + \|\sigma(Q,\om)
  \|^2) ds \} \; \|U\|^2,
\end{eqnarray}
where we used the Lemma \ref{lemma} on the LES model (\ref{LES}).
Also, by the Young and Poincar\'e inequalities \cite{Temam} we
have
\begin{eqnarray*}
\left|\beta \int_D (\bar{\psi}_x-\Psi_x)U \,dxdy \right| & \leq &
\frac12 |\beta| \left(\int_D (\bar{\psi}_x-\Psi_x)^2 \,dxdy +
\int_D U^2 \,dxdy\right)
 \\
& \leq &  \frac12 |\beta| \left(\frac{|D|}{\pi}\int_D U^2 \,dxdy +
\int_D U^2 \,dxdy \right),
\end{eqnarray*}
that is
\begin{equation}    \label{estimate4.5}
\left|\beta \int_D (\bar{\psi}_x-\Psi_x)U  \,dxdy \right| \leq
\frac12 |\beta| \left(\frac{|D|}{\pi} +1 \right) \|U\|^2.
\end{equation}
Moreover,
\begin{eqnarray}
 |\int_D [R(q)- \sigma(Q, \omega)]U| \leq \frac12 \|R(q)- \sigma(Q,
 \omega)\|^2 +\frac12 \|U\|^2.
\end{eqnarray}

Putting all these estimates into (\ref{Udot}), we obtain
\begin{eqnarray} \label{Udot2}
  \frac{d}{dt}\|U\|^2
  \leq & 2 &
  \left(- \a +\frac12 + \sqrt{\frac{2|D|}{\pi}} \;\right)\nonumber\\
&& \left(\|\bar{q}_0\|^2 e^{ 2|\a|T}
  +\frac1{r} e^{ 2|\a|T} \int_0^T \left[\|\bar{f}(s)\|^2 + \|\sigma(Q,\om)
  \|^2\right] ds \right) \;  \|U\|^2 \\
 && + \|R(q)- \sigma(Q,\omega)\|^2\nonumber.
\end{eqnarray}
Notice that $\a$ is defined in Lemma \ref{lemma} in terms of
physical parameters.
 By the Gronwall inequality \cite{Temam} and noting
$U(0)=0$, we obtain
\begin{equation}
\EX \|\bar{q}-Q\|^2 =\EX \|U\|^2 \leq C(\nu, r, q_0, T) \cdot\;
\EX \int_0^t \|R(q)- \sigma(Q,\omega)\|^2 dt, \;\;  0 \leq t \leq
T,
\end{equation}
where $C>0$ is a constant. This proves part (i) of Theorem 1.

To prove part (ii) of Theorem 1, denote $V= q -Q$, so that
$V= \D (\psi-\Psi)$.  Subtracting equation (\ref{QG}) from
(\ref{LES}) leads to
\begin{equation} \label{V}
 V_t =- J(\psi ,q) +  J(\Psi,Q) -\beta
 (\psi_x-\Psi_x)
  + \nu \Delta V - r V + (f -\bar{f})  - \sigma(Q, \omega).
\end{equation}
Similar to the approach in proving part (i) above, we estimate
\begin{eqnarray} \label{Vdot}
 \frac{d}{dt}\|V\|^2
 & \leq &
  2\left(- \a + 1 + \sqrt{\frac{2|D|}{\pi}} \right)\;
\left(\|q_0\|^2 e^{ 2|\a|T}
  +\frac1{r} e^{ 2|\a|T} \int_0^T \|f(s)\|^2 ds \right)
 \;  \|V\|^2  \nonumber \\
 && + \|f -\bar{f}\|^2 + \| \sigma(Q,\omega)\|^2.
\end{eqnarray}
By the Gronwall inequality again, we obtain
\begin{eqnarray}
\EX \|V\|^2 & \leq & C_1(\nu, r, q_0, T)\; \EX
\|q_0-\bar{q}_0\|^2  \nonumber \\
& + & C_2 (\nu, r, q_0, T) \cdot\;\EX \int_0^T \left[ \|f
-\bar{f}\|^2
 + \| \sigma(Q,\omega)\|^2 \right] dt, \;\;  0 \leq t \leq T,
\end{eqnarray}
where $C_1, C_2$ are positive constants. Due to the property of
$G_\delta$, both $\|q_0-\bar{q}_0\|$ and $ \|f
-\bar{f}\|$ go to zero as $\delta \to 0$. Together with the
condition (\ref{variance}), we finally see that $\EX \|V\|^2 =\EX
\|q-Q\|^2  \to 0$ as $\delta \to 0$, completing the proof of Theorem 1.

We view Theorem 1 as a starting point in our study of stochastic
parameterization of SGS stress related terms in the LES of geophysical
flows in the sense that it comments on the dependence of the LES
solution only on the variance of the stochastic parameterization (its
closeness to that of the actual SGS stress terms).  Thus, it may be
argued that the stochastic nature of the parameterization is not
central to the results of this paper. However, to the extent that
variance is one of the most fundamental characteristics of a
stochastic process, understanding the dependence of the LES solution
on it is important. What is now desirable is further characterization
of the dependence of the LES solution on various other aspects of the
stochastic parameterization such as its temporal and spatial
correlation structure and its probability distribution function. These
are subjects of ongoing research and we hope to report on them in the
future.


{\bf Acknowledgement.}
  This work was initiated while both authors were co-organizing
the special session on \emph{Uncertainty, Random Dynamical Systems
and Stochastic Modeling in Geophysics}, European Geosciences
Union, General Assembly,   Vienna, Austria,  April 2005.
   We thank Traian Iliescu for helpful discussions.


\bibliographystyle{amsplain}

\end{document}